\documentclass[11pt]{amsart}
\usepackage{vmargin}
\usepackage{harvard}
\usepackage{psfig}
\usepackage{graphics}

\newtheorem{theorem}{Theorem}[section]

\newtheorem{lemma}[theorem]{Lemma}

\newtheorem{proposition}[theorem]{Proposition}

\newtheorem{corollary}[theorem]{Corollary}

\theoremstyle{definition}
\newtheorem{definition}[theorem]{Definition}
\newtheorem{example}[theorem]{Example}

\theoremstyle{remark}

\newtheorem*{unremark}{Remark}

\makeatletter
\def\romenumi{%
  \def\theenumi{\roman{enumi}}%
  \def\p@enumi{\theenumi}%
  \def\labelenumi{(\@roman\c@enumi)}}
\makeatother

\newcommand{\noproof}{\hfill\qedsymbol}

\DeclareMathOperator{\re}{Re}
\DeclareMathOperator{\im}{Im}
\DeclareMathOperator{\Arg}{Arg}

\DeclareMathOperator{\hess}{\mathcal{H}}

\DeclareMathOperator{\grad}{\nabla}
\DeclareMathOperator{\sign}{sgn}

\numberwithin{equation}{section}

\setpapersize{USletter}
\setmarginsrb{20mm}{20mm}{20mm}{20mm}{12pt}{10mm}{12pt}{10mm}

\def\sing{\mathcal{V}}
\def\torus{T}
\def\disk{D}
\def\numer{{G}}
\def\denom{{H}}

\def\direc{{\bf dir}}
\def\dom{{\mathcal D}}
\def\logdom{{\rm log} \dom}

\def\nbd{{\mathcal N}}
\def\one{{\bf 1}}
\def\zero{{\bf 0}}
\def\GF{{\mbox{GF-sequence}}}
\def\cio{C^\infty_0}
\def\finv{\eta}
\def\X{\Xi}

\def\b#1{\mathbf{#1}}

\def\bt{{\tilde{b}}}
\def\gt{{\tilde{g}}}
\def\ft{{\tilde{f}}}
\def\psit{{\tilde{\psi}}}
\def\R{{\mathbb{R}}}
\def\CC{{\mathbb{C}}}
\def\RP{{\mathbb{RP}}}
\def\N{{\mathbb{N}}}

\def\A{{\mathcal A}}

\def\C{{\mathcal C}}

\def\ee{\varepsilon}

\begin{document}

\title[Asymptotics of multivariate sequences I]{Asymptotics of multivariate sequences, part I: smooth points of the singular variety}
\author{Robin Pemantle}
\address{Department of Mathematics, Ohio State University, Columbus OH 43210}
\email{pemantle@math.ohio-state.edu}
\thanks{Research supported in  part by NSF grant DMS-9803249}
\author{Mark C. Wilson}
\address{Department of Mathematical Sciences, University of Montana,
 Missoula, MT 59812}
\email{wilsonm@member.ams.org}



\begin{abstract}Given a multivariate generating function
$F(z_1 , \ldots , z_d) = \sum a_{r_1 , \ldots , r_d} z_1^{r_1}
\cdots z_d^{r_d}$, we determine asymptotics for the coefficients.
Our approach is to use Cauchy's integral formula near singular
points of $F$, resulting in a tractable oscillating integral. This
paper treats the case where the singular point of $F$ is a smooth
point of a surface of poles. Companion papers will treat singular
points of $F$ where the local geometry is more complicated, and
for which other methods of analysis are not known.
\end{abstract}

\maketitle

\section{Introduction}

The generating function $F(z) := \sum_{r=0}^\infty a_r z^r$ for
the sequence $a_0 , a_1 , a_2 , \ldots$ is one of the most useful
constructions in combinatorics.  If the function $F$ has a simple
description, it is usually not too hard to obtain $F$ as a formal
power series once one understands a recursive or combinatorial
description of the numbers $\{ a_r \}$.  One may then analyze the
analytic properties of $F$ in order to obtain asymptotic
information about the sequence $\{ a_r \}$.  While still part art
and part science, this latter analytic step has become quite
systematized. \citeasnoun{stanley;enumerative-combinatorics1} in
his introduction to enumerative combinatorics gives the example of
the function $F(z) = \exp (z + \frac{z^2}{2})$, from which he says
``it is routine (for someone sufficiently versed in complex
variable theory) to obtain the asymptotic formula $a_r = 2^{-1/2}
r^{r/2} e^{-r/2 + \sqrt{r} - 1/4}$.'' Routine, in this case, means
a single application of the saddle point method.  When $F$ has
singularities in the complex plane, the analysis is often more
direct: the location of the singularities and the behavior of $F$
near these determine almost algorithmically the asymptotic
behavior of the sequence $\{ a_r
\}$.  For those not sufficiently versed in complex variable
theory, two useful sources are \citeasnoun{henrici;complex-analysis2} and
\citeasnoun{odlyzko;asymptotic-methods}.  The transfer theorems of
\citeasnoun{flajolet-odlyzko;singularity-analysis}
encapsulate much of this knowledge in a very useful way; see also
\citeasnoun{wilf;GFology} for an elementary introduction.

When the sequence $a_0 , a_1 , a_2 , \ldots$ is replaced by a
multidimensional array $\{ a_{r_1 , \ldots , r_d} \}$, things
become much more hit and miss.  Let us use boldface to denote
vectors in $\CC^d$ or $\N^d$, and use multi-index notation, so
that $a_\b{r}$ denotes the multi-index $a_{r_1 , \ldots , r_d}$ and
$\b{z}^\b{r}$ denotes the product $z_1^{r_1} \cdots z_d^{r_d}$ which
we will sometimes write in expanded form for clarity.  The
generating function $F : \CC^d \to \CC$ is defined analogously to
the one-dimensional generating function by
$$F(\b{z}) = \sum_{\b{r} \in \N^d} a_\b{r} \b{z}^\b{r} .$$
Surprisingly, techniques for extracting asymptotics of $\{ a_\b{r}
\}$ from the analytic properties of $F$ were, until recently,
almost entirely missing.  In a survey of asymptotic methods,
\citeasnoun{bender;asymptotic-methods} says:
\begin{quote}
Practically nothing is known about asymptotics for recursions in
two variables even when a generating function is available.
Techniques for obtaining asymptotics from bivariate generating
functions would be quite useful.
\end{quote}
In the intervening 25 years, some results have appeared,
addressing chiefly the case where the array $\{ a_\b{r} \}$ obeys a
central limit theorem.  Common to all of these is the following
method.  Treat $\{ a_\b{r} \}$ as a sequence of
$(d-1)$-dimensional arrays indexed by $r_d$; show that the $n^{th}$
$(d-1)$-dimensional generating function is roughly the $n^{th}$
power of a given function; use this approximation to invert the
characteristic function and obtain a Central Limit Theorem.  We
refer to these methods as $\GF$ methods.  The other body of work
on multivariate sequences, which we will call the diagonal method,
is based on algebraic extraction of the diagonal, as found in
 \citeasnoun{hautus-klarner;diagonal} (  
see also \citeasnoun{furstenberg;diagonal} and later \citeasnoun{lipshitz;diagonal} for an algebraic description of the scope of this method; variants are described in 
\citeasnoun{stanley;enumerative-combinatorics2} and \citeasnoun{pippenger;equicolourable}).

The most fundamental $\GF$ result is probably
\citeasnoun{bender-richmond;limit-multivariate2},
with extensions appearing in later work of the same
authors.  \citeasnoun{flajolet-sedgewick;multivariate} present a version of the
same idea which holds in much greater generality.
\citeasnoun{gao-richmond;limit-multivariate4} go beyond the central limit case,
using the transfer theorems of
\citeasnoun{flajolet-odlyzko;singularity-analysis} to handle
functions that are products of powers with powers of logs. Recent
work of Bender and Richmond
\cite{bender-richmond;admissible,bender-richmond;products}
extends the applicability of the central limit results to many
problems of combinatorial interest; see also
\cite{hwang;stirling,hwang;convergence-rates}, where
more precise asymptotics are given, and
\citeasnoun{hwang;local-limit-deviations}, which extends some results
to the combinatorial schemes of
\citeasnoun{flajolet-soria;gaussian-exponential}. This does not
exhaust the recent work on the problem of multivariable
coefficient extraction, but does circumscribe it.

The present paper, together with forthcoming companion papers,
takes aim at a large class of multivariable coefficient extraction
problems, for which a fair amount of information can be read off
in a systematic way.  An ultimate goal (not our only goal) is to
systematize the extraction of multivariate asymptotics
sufficiently that it may be automated, say in Maple.  Everything
we do, we do with complex contour integration.  In this regard,
our methods are most similar to those of
\citeasnoun{bertozzi-mckenna;queueing}, who, as we do, provide a general
framework for harnessing the multivariable theory of residues for
exact and series computation of coefficients. A more detailed
description of our method will be given in
Section~\ref{ss:results}, but here is an outline.

\begin{quote}
 ~~

(1) Use the multidimensional Cauchy integral formula to represent
$a_{\b{r}}$ as an integral over a $d$-dimensional torus inside
$\CC^d$.

(2) Expand the surface of integration across a point $\b{z}$ where
$F$ is singular, and use the residue theorem to represent
$a_\b{r}$ as a $(d-1)$-dimensional integral of one-variable residues.
The choice of $\b{z}$ determines the directions in which asymptotics
may be computed.

(3) Put this in the form of an integral $\int \exp (\lambda
f(\b{z})) \psi (\b{z}) \, d\b{z}$ for which the large-$\lambda$
asymptotics can be read off from the theory of oscillating
integrals.
\end{quote}

\begin{figure}
\centerline{\psfig{figure=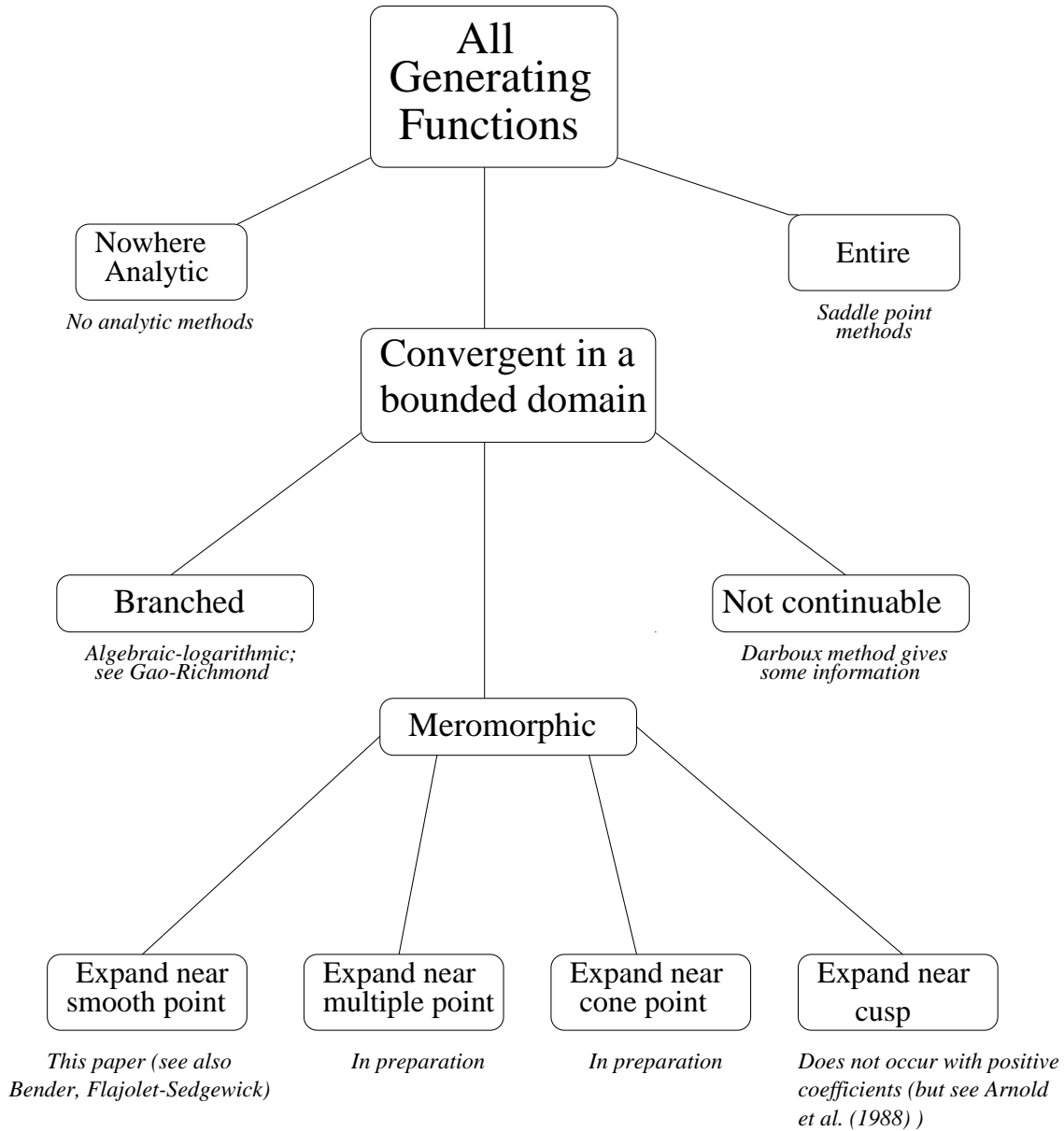,width=6in}}
\vspace{10pt}
\caption{Classification of generating functions}
\label{fig:taxonomy}
\end{figure}

In the rest of this introductory section, we describe the scope of
our methods. Figure~\ref{fig:taxonomy} depicts a classification of generating functions and illustrates the remainder of this paragraph.  If a
formal power series is nowhere convergent, analytic methods are
useless. Among those power series converging in some neighborhood
of the origin, there are three possibilities: a function may be
entire, may have singularities around which analytic continuations
exist, or it may be defined only on some bounded subset of
$\CC^d$.  Our methods are tailored to the second class.  The third
class, although in some sense generic, seldom arises in any
problem for which the generating function may be effectively
described.  Incomplete asymptotic information is available via
Darboux' method; details of this method in the univariate case are
given in \citeasnoun{henrici;complex-analysis2} and
\citeasnoun{odlyzko;asymptotic-methods}.  The first class can
and does arise frequently. Our methods are simply not equipped to
handle entire functions, and systematizing the asymptotic analysis
of coefficients of entire generating functions remains an
important open problem.

For the remainder of this paper, we will assume that the formal
power series $F$ converges in a neighborhood of the origin and may
be analytically continued everywhere except a set $\sing$ of complex
dimension $d-1$ which we call the {\em singular variety}. The
point $\b{z}$ in step~2 is an element of $\sing$, and the behavior of
$\sing$ near $\b{z}$ greatly affects the subsequent analysis in step~3.
This paper addresses the case where $\b{z}$ is a smooth point of
$\sing$ at which $F$ has a pole.  The forthcoming companion papers
will address cases where $\b{z}$ is a multiple point
or a cone point.  We do not know whether cases
where $\b{z}$ is a cusp of $\sing$ arise, but if so, the subsequent
analysis has mostly been carried out in the work of
\citeasnoun{arnold-varchenko;asymptotics-integrals}.

The chief purpose of this study is to give a solution to the
problem of asymptotic evaluation of coefficients that is as
general as possible.  An important part of this is re-derivation
in a general setting of results obtainable via $\GF$ or {\em ad
hoc} methods. We show in Section~\ref{ss:classify} how unifying
these results allows us to show that our method successfully finds
asymptotics for every function in a certain large class.  Familiar
examples from this class include: lattice path counting, various
known generating functions for polyominos and stacked balls,
enumeration of Catalan trees by number of components or
surjections by image cardinality (see
\citeasnoun{flajolet-sedgewick;multivariate},
stopping times for certain random walks
(see \citeasnoun{larsen-lyons;coalescing}), as well as the examples
given in the $\GF$ papers of
\citeasnoun{bender;central-local} and
\citeasnoun{bender-richmond;limit-multivariate2}: ordered set partitions
enumerated by number of blocks, permutations enumerated by rises,
and Tutte polynomials of recursive sequences of graphs.

Nevertheless, our pursuit of this problem was also motivated by
some specific applications which we mention briefly now and
discuss more thoroughly later.  These are cases where known
methods do not suffice to obtain complete asymptotic information.
There is a class of tiling enumeration problems for which an
explicit three variable rational generating function may be
obtained.  This class includes the Aztec Diamond domino tilings of
\citeasnoun{cohn-elkies-propp;aztec}. Asymptotics in the so-called {\em
region of fixation} are obtained from analysis of the smooth
points of $\sing$ (Theorem~\ref{th:higher d} below), while
asymptotics in the region of positive entropy are derived from
analysis of the cone point.  \citeasnoun{cohn-pemantle;fixation}  applies a
cone point analysis to a tiling enumeration problem for which the
only previous results are some pictures via simulation
(http://www.math.harvard.edu/\~{}cohn/picture.gif) . Another
motivation has been to solve the general multivariable linear
recursion. Depending on whether one allows forward recursion in
some of the variables, one obtains either rational or algebraic
generating functions.  The general rational function may have any
of the types of singularities mentioned above: smooth points,
nodes, cones, cusps, branchpoints, etc. Even the simple rational
generating function $1 / (3 - 3z - w + z^2)$ of
Example~\ref{ex:cubic} requires two separate analyses in order to
get asymptotics in all directions. We will see that
Theorem~\ref{th:simple} gives asymptotics in one region, while
Theorem~\ref{th:2d} is required for other directions.

Asymptotics derived near smooth pole points nearly always exhibit
central limit behavior.  Smooth pole points are the topic of this
first paper, and are exactly the case to which existing methods
may apply.  While one function of this paper is to lay foundations
for the cases in which the singularity is more complicated, there
are several ways in which it improves upon available analyses of
the smooth case.

First, most of the existing results assume that the singular point
$\b{z} \in \sing$ has positive real coordinates, and that it is
strictly minimal in a sense defined in the next section.  This
assumption often holds when the coefficients $\{ a_\b{r} \}$ are
nonnegative reals, though it will fail if, for example, there is
any periodicity.  The assumption always fails when the
coefficients $\{ a_\b{r} \}$ have mixed signs, as is the case for
example with the generating functions $(1 - zw)/(1 - 2zw +w^2)$
and $1 / (1 - 2zw + w^2)$ for the Chebyshev polynomials of the
first and second kinds \cite[page 50]{comtet;advanced}. $\GF$ methods may
be adapted to some of these situations. Indeed, the presentation
of these methods by \citeasnoun[Theorem~9.7]{flajolet-sedgewick;multivariate}
accomplishes this adaptation in great generality. But certainly
there are cases such as the rational generating function
$1 / (1 - z - w + \beta zw)$, where the points $\b{z}$ with given
moduli form a continuum and standard $\GF$ methods are not
sufficient.

Second, our methods obtain automatically a full asymptotic
expansion of $a_{r_1 , \ldots , r_d}$ in decreasing powers of the
indices $r_j$.  This is certainly not inherent in the existing
results, whose relatively short proofs involve inversion of the
characteristic function (see however
\citeasnoun{hwang;stirling} and \citeasnoun {hwang;central-limit-deviations} for
something in this direction). The expansion to $n$ terms is
completely effective in terms of the first $n$ partial derivatives
of $1/F$ at $\b{z}$, as is the error bound.

Third, these results explicitly cover the case where the pole at
$\b{z}$ has order greater than 1.  The behavior in this case is not
according to the central limit theorem.  The only existing work
addressing this case is \citeasnoun{gao-richmond;limit-multivariate4},
and they require
nonnegativity assumptions, as mentioned above.  In the case where
$F = G^k$ is an exact power, one could attempt first to solve the
problem for $G$ and then to take the $k$-fold convolution.  This
is much harder than the present approach, as may be seen by the
rather involved computation in \citeasnoun{cohn-elkies-propp;aztec}.

Fourth, the potential for increasing the scope to new applications
seems greater for contour methods than for $\GF$ methods.  The
contour method reduces the asymptotic problem to the problem of an
oscillating integral near a singularity, which can almost
certainly be done.  By contrast, the $\GF$ method requires first
an understanding of the sequence of $(d-1)$-dimensional generating
functions arising from the given $d$-dimensional generating
function, and then another result in order to transfer this
information to asymptotics of the coefficients $a_\b{r}$.

Fifth, although our results in the case of smooth pole points are
often similar to those obtained by $\GF$ methods, our hypotheses
are quite different.  In Section~\ref{ss:classify} we show how our
hypotheses may be universally established for functions that
generate nonnegative values and are meromorphic through their
domain of convergence.

Finally, we compare our method to recent results from the diagonal
method.  It is known \cite{lipshitz;diagonal} that the diagonal
sequence
$a_{n,n,\ldots ,n}$ of a multivariate sequence with rational
generating function has a generating function satisfying a linear
differential equation over rational functions.  Much is known
about how to compute this equation (see for example
\citeasnoun{chyzak-salvy;symbolic-integration}).
If one wants asymptotics on the diagonal, or in
any direction where the coordinate ratios are rational numbers
with small denominators, then these methods give results that are
in theory at least as good as ours.  The method, however, is
inherently non-uniform in the direction, so there is no hope of
extending it to larger sets of directions, which is what we
accomplish in the present work.

The remainder of the paper is organized as follows.  In the next
section we set forth notation and define the terms necessary to
state the main results of the paper.  The main results are stated
in Section~\ref{ss:results}, and examples are given.  The next
section contains a proof of these results, modulo the computation
of some oscillating integrals.  This computation is carried out in
Section~\ref{ss:oscillate}.  Section~\ref{ss:classify} outlines
some details of taxonomy and discusses universality of the method
of complex contour integration.  The final section states some
open problems.

\section{Notation and Preliminaries}

The main results of this paper give asymptotics valid under
certain geometric assumptions on $\sing$ and computable from some
quantities that are in turn effectively computable from the
generating function $F$. Thus in addition to setting out basic
notation, we need to define some terms related to the geometry of
$\sing$ and some quantities associated with $F$.

\subsection{Notation}

Throughout the paper, $F = \sum a_\b{r} \b{z}^\b{r}$ will denote a
function on $\CC^d$ analytic in a neighborhood of the origin. The
(open) domain of convergence of the power series will be denoted
$\dom$. For $\b{z} \in \CC^d$, let $\torus (\b{z})$ denote the torus
consisting of points $\b{w}$ with $|w_j| = |z_j|$ for $1 \leq j \leq
d$ and let $\disk (\b{z})$ denote the closed polydisk of points
$\b{w}$ with $|w_j| \leq |z_j|$ for $1 \leq j \leq d$.  Recall (see
\citeasnoun{hormander;intro}) that the domain $\dom$ is a union of tori
$\torus (\b{z})$ and is logarithmically convex, that is, the set
$$\logdom := \{ \b{x} \in \R^d : (e^{x_1} , \ldots , e^{x_d})
   \in \dom \}$$
is a convex subset of $\R^d$ and is an order ideal, that is, it is
closed under $\leq$ in the coordinatewise partial order.

We assume throughout that $F = \numer / \denom$, where both
$\numer$ and $\denom$ are analytic in a neighborhood of
$\disk (\b{z})$ for some point $\b{z}$.  In particular, every
meromorphic function satisfies this condition
\footnote{The greater
generality allows us to cover examples such as the generating
function for self-avoiding random walks
\cite{chayes-chayes;self-avoiding} or percolation paths in the
subcritical regime \cite{campanino-chayes-chayes;percolation}.  In
these cases, all the work is in showing the function is
meromorphic in a neighborhood of
$\disk (\b{z})$. Without further knowledge, the authors then
conclude central limit behavior.}.  The set where
$\denom$ vanishes will be denoted $\sing$.  Many of our examples will
be in dimension 2, in which case we will often use $z$ and
$w$ in place of $z_1$ and $z_2$, use $(z,w)$ in place of
$\b{z}$, and use $(r,s)$ in place of $(r_1 , r_2)$.  We sometimes
need to treat $\CC^d$ as $\CC^{d-1} \times \CC$ (although symmetry of
the coordinates is preserved most of the time).  Accordingly, when
the dimension is greater than~2, we use $\widehat{\b{z}}$ to denote
$(z_1 , \ldots , z_{d-1})$.  Partial derivatives will be denoted
$H_1$ for $\frac{\partial H}{\partial z_1}$ and so forth;
in dimension~2 we will also use $H_z$ and $H_w$.

As is usual for asymptotic analyses, we let $f \sim g$ denote
$f / g \to 1$, with the limit taken at infinity unless otherwise
specified.  The function $f$ is said to be {\em rapidly
decreasing} if $f(x) = O(x^{-N})$ for every $N$, and is said to be
{\em exponentially decreasing} if $f(x) = O(e^{-cx})$ for some $c
> 0$. We also use the symbol ``$\sim$'' to denote asymptotic
expansion. Thus
$$f \sim \sum b_n g_n$$
is normally taken to mean that $f - \sum_{n=0}^N b_n g_n = o(b_N
g_N)$, where $b_n\in \mathbb{C}$ and $\{ g_n \}$ is a fixed sequence of functions such that
$g_{n+1} = o(g_n)$ for each $n$. We broaden this to allow $b_n = 0$
when $n \neq 0$, so that the remainder term need only be $o(g_n)$
and not $o(b_n g_n)$.  In particular, if
$$f(x) \sim g(x) \cdot \sum_{n = 0}^\infty c_n x^{-n}$$
with $c_0 = 1$, then we say we have obtained a full asymptotic
expansion for $f$ in decreasing powers of $x$ with leading term
$g$.

\subsection{Geometry of $\sing$} \label{ss:geom}

As in the one-dimensional case, the points of $\sing$ nearest the
origin are the most important.  Accordingly we define a point $\b{z}
\in \sing$ to be {\em minimal} if $\sing \cap \disk (\b{z}) \subseteq
\torus (\b{z})$; we say that $\b{z}$ is {\em locally minimal}
if the analogous relation holds with $\sing$ replaced by a neighborhood of $\b{z}$ in $\sing$. Divide the minimal
points of $\sing$ into three types. Say that
$\b{z}$ is {\em strictly minimal}, {\em finitely minimal} or {\em
toral}, according to whether the cardinality of
$\sing \cap \disk (\b{z})$ is 1, finite, or infinite.  When infinite, the
intersection must be uncountable.  If $\b{z}$ is a minimal point of
$\sing$ then the interior of $\disk (\b{z})$ is contained in $\dom$,
so the assumption that $\numer$ and $\denom$ are analytic on a
neighborhood of $\disk (\b{z})$ is just a little stronger than what
is true automatically.

A {\em simple pole} of $F$ is a point $\b{z} \in \sing$ where $\denom$
vanishes to order 1.  Equivalently, the gradient $\grad \denom$
does not vanish.  Let $\b{z}$ be a simple pole of $F$ and assume for
specificity that $\denom_d$ is nonzero at $\b{z}$.   By the implicit
function theorem, there is a neighborhood of $\b{z}$ where
$\sing$ may be parametrized by $z_d = g(z_1 , \ldots , z_{d-1})$ for
some analytic function $g$.  We will always use $g$ to denote this
parametrization.

We will see later (in the proof of Theorem~\ref{th:classify}) that
under some hypotheses on $F$, minimal points of $\sing$ are always
found in the positive real orthant.  A relation true in complete
generality is the following.

\begin{lemma} \label{lem:dir}
Let $\b{z}$ be a simple pole of $F$ and suppose that $z_d \denom_d$
does not vanish there.  If $\b{z}$ is locally minimal then for all
$j < d$, the quantity $z_j \denom_j / (z_d \denom_d)$ is real and
nonnegative.
\end{lemma}

\begin{proof} Given $\theta$ and $j$, let $\b{z}^{(\theta)}$
be the result of varying $\b{z}$ by multiplying the $j^{th}$
coordinate by $e^{i \theta}$ and adjusting the last coordinate so
as to remain on $\sing$ (that is, \hfill \\ $z^{(\theta)}_d = g(z_1 ,
\ldots , z_{j-1} , z_j e^{i \theta} , z_{j+1} , \ldots ,
z_{d-1})$). Differentiating the relation $\denom (\b{z}^{(\theta)})
= 0$ implicitly with respect to
$\theta$ at 0 yields
\begin{equation} \label{eq:minimal}
i z_j \denom_j + \denom_d {d z^{(\theta)}_d \over d \theta} = 0 .
\end{equation}
By minimality of $\b{z}$, we know that the modulus of
$z^{(\theta)}_d$ has a minimum at $\theta = 0$, hence
$(d z^{(\theta)}_d / d \theta) / z_d$ is purely imaginary.
Plugging this into~(\ref{eq:minimal}) proves that $z_j \denom_j /
(z_d \denom_d)$ is real.  If $z_j \denom_j / (z_d \denom_d) =
-\beta < 0$ then $\sing$ has a tangent vector at $\b{z}$ in the
direction $- z_j e_j - \beta z_d e_d$, where $e_j$ is the
$j^{th}$ coordinate vector.  This contradicts minimality.
Hence $z_j \denom_j / (z_d \denom_d) \geq 0$.   \end{proof}

\begin{definition}
Define $\direc (\b{z})$ to be the equivalence class of (complex)
scalar multiples of the vector $(z_1 \denom_1 , \ldots , z_d
\denom_d)$, defined whenever $z_j \denom_j$ does not vanish for
all $j$. By the previous lemma, when $\b{z}$ is a minimal pole of
$F$ with nonzero coordinates, $\direc (\b{z})$ is a well defined
element of
$\RP^{d-1}$.
\end{definition}
The importance of $\direc$ is that analysis of $F$ near $\b{z}$
yields asymptotic information about $a_\b{r}$ with $\b{r}
\in \direc (\b{z})$.  The function $\direc$ appears in $\GF$
method literature as ${\bf m}$.  When $\b{z} \in \partial \dom$ is
on the boundary of the domain of convergence, $\direc (\b{z})$ is
the normal to the support hyperplane of the convex set $\logdom$
at the point $(\log |z_1| , \ldots , \log |z_d|)$.

We now define a few more quantities associated with $F$ and $g$.
Again, we will reserve the names of these functions, so as not to
burden the notation with subscripts and arguments.  If $\b{z}$ is a
simple pole of $F$ with $z_d \denom_d$ not vanishing there, define
a function $\psi$ on a neighborhood of $\widehat{\b{z}}$ by
\begin{equation} \label{eq:psi}
\psi (\widehat{\b{w}}) = - \lim_{w \rightarrow g(\widehat{\b{w}})}
   (w - g(\widehat{\b{w}})) {F (\widehat{\b{w}} , w) \over w} .
\end{equation}
Suppose now that $\widehat{\b{w}} \in \torus (\widehat{\b{z}})$ and
write
$w_j = z_j e^{i \theta_j}$.  For fixed $\b{r}$ with $r_d \neq 0$,
define a function $f$ on a neighborhood of $\widehat{\b{z}}$ in
$\torus (\widehat{\b{z}})$ by
\begin{equation} \label{eq:f}
f (\widehat{\b{w}}) = \log \left ( {g (\widehat{\b{w}}) \over
   g(\widehat{\b{z}})} \right ) + i \sum_{j=1}^{d-1} {r_j \over r_d}
   \theta_j .
\end{equation}

We will be parametrizing integrals over $\torus (\widehat{\b{z}})$
by ${\bf \theta}$, so we will want the above function expressed in
terms of $\widehat{\bf \theta}$.  We therefore compose with the
map $M$ taking $\widehat{\bf \theta}$ to $\widehat{\b{w}}$ defined
by
$M(\theta_1 , \ldots , \theta_{d-1}) = (z_1 e^{i \theta_1} ,
\ldots , z_{d-1} e^{i \theta_{d-1}})$, and define the functions
$\gt: = g \circ M , \ft := f \circ M , \psit := \psi \circ M$.

Although it is not obvious yet, $\ft$ will always vanish at
$\zero$  to at least two orders (Lemma~\ref{lem:stat} below), and
the hypothesis $Q \neq 0$ in Theorem~\ref{th:simple} is equivalent
to $\ft$ having nonvanishing quadratic term.  For ease of
reference, Table~\ref{table1} summarizes the foregoing definitions,
stratified by how many times the given data $\numer$ and $\denom$
have been manipulated.


\begin{table}
\label{table1}
\caption{Reserved notation in remainder of this article}
\begin{tabbing}
Given \= information: \\
\> the function $F$ in the form $\numer / \denom$ \\[3ex]
First level: \\
\> $g$ parametrizes the zero set, $\sing$ of $\denom$ \\
\> $\direc(\b{z})$ is the coordinatewise product $(\grad \denom)\cdot(\b{z})$
   in projective space \\[3ex]
Second level: \\
\> $\psi$ is the residue in $z_d$ of $F / z_d$ at points $(\widehat{\b{z}} ,
g(\widehat{\b{z}}))$ \\
\> $f$ is $\log g$, plus a term linear in $\log z_j$ and depending
   on $\b{r}$. \\[3ex]
Third level: \\
\> $\psit , \gt$ and $\ft$ are $\psi, g$ and $f$ expressed in
   terms of ${\bf \theta}$
\end{tabbing}
 ~~\\[2ex]
\end{table}

\section{Statement of results, with examples} \label{ss:results}

Before going on, we pause to state a prototype of our results in
the simplest possible setting, namely where the number of
variables is~2, the functions $\numer$ and $g$ are as
nondegenerate as possible, and only the leading term asymptotic is
given.  The proof is in Section~\ref{ss:proofs}.

\begin{theorem} \label{th:simple}
Let $F = \numer / \denom$ be a meromorphic function of two
variables, not singular at the origin.  Define
$$Q (z,w) := - w^2 \denom_w^2 z \denom_z - w \denom_w z^2 \denom_z^2 - w^2 z^2
   \left ( \denom_w^2 \denom_{zz} + \denom_z^2 \denom_{ww}
   - 2 \denom_z \denom_w \denom_{zw} \right ) .$$
Then
$$a_{r,s} \sim {\numer (z,w) \over \sqrt{2 \pi}} z^{-r} w^{-s}
   \sqrt{- w \denom_w \over s Q}$$
uniformly as $(z,w)$ varies over a compact set of strictly
minimal, simple poles of $F$ on which $Q$ and $\numer$ are
nonvanishing, and $(r,s) \in \direc (z,w)$.
\end{theorem}

\noindent{\em Remarks:} Usually the expression in the radical will
be positive real, as will the coefficients $a_{rs}$.  The result
is true in general, though, as long as the square root is taken to
be $-w \denom_w$ times the principal root of $Q / (- w
\denom_w^3)$. Also note that when $(r,s) \in \direc (z,w)$ then
the expression
$w \denom_w / s$ is coordinate-invariant, that is, equal to $z \denom_z / r$.
Thus the given expression for $a_{r,s}$ has the expected symmetry.

\begin{example}[Lattice paths] \label{ex:lattice} ~~ \end{example}

Let $a_{r,s}$ be the number of nearest-neighbor paths from the
origin to $(r,s)$ moving only north, east and northeast; these are
sometimes called {\em Delannoy numbers}~\citeasnoun[page
185]{stanley;enumerative-combinatorics2}. The generating function
is $F(z,w) = 1 / (1 - z - w - zw)$. The zero set $\sing$ of
$\denom = 1 - z - w - zw$ is given by
$w = (1-z)/(1+z)$, and the minimal points of $\sing$ are those where
$w \in [0,1]$.  With the help of relations that hold when $\b{z} \in \sing$
we may compute as follows.
\begin{eqnarray*}
\denom_z & = & - 1 - w \\
   -z \denom_z & = & 1 - w \\
   Q & = & (1-z) (1-w) (1-zw) \\
   {z \denom_z \over w \denom_w} & = & {1-w \over 1-z} =
   {1 - w^2 \over 2w}
\end{eqnarray*}
with $\denom_w$ and $-w \denom_w$ given by reversing $z$ and $w$.
As $z$ varies over $[\ee , 1-\ee]$, the functions $Q$ and $\numer
:= 1$ do not vanish.  The minimal pair $(z,w)$ that solves $(r,s)
\in \direc (z,w)$ is given by $z = (\sqrt{r^2 + s^2} - s)/r$ and
$w = (\sqrt{r^2 + s^2} - r)/s$.  Theorem~\ref{th:simple} then
gives
\begin{eqnarray*}
a_{rs} & \sim & \left ({\sqrt{r^2 + s^2} - s \over r} \right)^{-r}
   \left ( {\sqrt{r^2 + s^2} - r \over s} \right )^{-s}
   \sqrt{1 \over 2 \pi} \sqrt{{1-z \over s} {1 \over 1- zw}}
   \\[1ex]
& = & \left ({\sqrt{r^2 + s^2} - s \over r} \right)^{-r}
   \left ( {\sqrt{r^2 + s^2} - r \over s} \right )^{-s}
   \sqrt{1 \over 2 \pi} \sqrt{r s \over (r+s-\sqrt{r^2 + s^2})^2
   \sqrt{r^2 + s^2}} \, ,
\end{eqnarray*}
uniformly when $r/s$ and $s/r$ remain bounded.  In particular,
when $r=s=n$, this gives the following formula for the $n^{th}$
diagonal coefficient (which may alternatively be obtained by
computing the diagonal generating function $(1 - 6s+ s^2)^{-1/2}$
according to the method given in
\citeasnoun[Section~6.3]{stanley;enumerative-combinatorics2}:
$$(\sqrt{2} - 1)^{-2n} \sqrt{1\over 2 \pi} {2^{-1/4} \over
  2 - \sqrt{2}} \, .$$

The computations in Theorem~\ref{th:simple} in terms of the values
and derivatives of $\numer$ and $\denom$ are explicit.  As we
state more general theorems, it becomes cumbersome and in fact
obfuscating to give formulae for the expansion coefficients
directly in terms of derivatives of $\numer$ and $\denom$.  This
is one reason we have already introduced the functions in
Table~\ref{table1}. It should be emphasized, however, that while we use
higher level quantities in the statements of subsequent theorems,
each expansion coefficient can be computed from finitely many
derivatives of $\numer$ and $\denom$.  We begin with a relatively
explicit computation for the general two-variable case.

For $k$ at least 2, we define constants
\begin{eqnarray}
A_+ (k,l) & := & {1 \over k} \Gamma \left ({l+1 \over k} \right )
 \label{eq:plus} \\[2ex]
A (k,l) & := & {1 \over k} \Gamma \left ( {l+1 \over k} \right )
   \left ( 1 + e^{\sign\Arg (c_k) i \pi (l - {l+1 \over k} )} \right ) \,
   \mbox{ if } k \mbox{ is odd}, \label{eq:k odd} \\[2ex]
A (k,l) & := & {2 \over k} \Gamma \left ( {l+1 \over k} \right )
   \mbox{ if } k,l \mbox{ are even},\label{eq:k even} \\[2ex]
A (k,l) & := & 0 \mbox{ if } k \mbox{ is even and } l \mbox{ is
   odd.} \nonumber
\end{eqnarray}
Let
$$y (x) = f(x)^{1/k} = c_k^{1/k} x \left ( 1 + {f(x) - c_k x^k
   \over c_k x^k} \right )^{1/k} \, ,$$
where $c_k$ is the first nonvanishing Taylor coefficient of $f (x)
= \sum_{j=k}^\infty c_j x^j$ and the argument of
$c_k^{1/k}$ is taken between $-\pi/(2k)$ and
$\pi / (2k)$.  Let $\finv$ denote the inverse function to $y$
and let $\{ b_j \}$ be the Taylor coefficients of
$(\psit \circ \finv) \cdot \finv'$. Clearly each $\{ b_j \}$
is determined by finitely many partial derivatives of $\numer$ and
$\denom$, and the index $l_0$ of the first nonvanishing $b_l$ is
the same as the order of vanishing of $\psit$ at 0.  The
coefficients $b_l$ are easily computed from the coefficients
$\bt_j := \psit^{(j)} (0) / j!$ and $c_j := \ft^{(j)} (0) / j!$;
in particular, if $\ft \sim c_k x^k$ near 0 then
\begin{equation} \label{eq:b}
b_{l_0} = \bt_{l_0} c_k^{-1/k} \, .
\end{equation}

\begin{theorem} \label{th:2d}
Let $F = \numer / \denom = \sum a_{rs} z^r w^s$ have a strictly
minimal, simple pole at $(z,w)$. Let $k$ be the order of vanishing
of $\ft$ at 0.  Let $l_0$ be the order to which $\numer$ vanishes
near $(z,w)$ on $\sing$, that is, the largest $l$ such that $G(z',w')
= O(|z-z'|^l + |w-w'|^l)$ as $(z' , w') \to (z,w)$ in $\sing$.  Then
there is a full asymptotic expansion
\begin{equation}\label{eq:full 2d}
a_{r,s} \sim {1 \over 2 \pi} z^{-r} w^{-s} \sum_{l \geq l_0}
   \A (k,l) b_l s^{-(l+1) / k} \, ,
\end{equation}
where $\A (k,l)$ denotes $A (k,l)$ if $\im \{ c_k \} \geq 0$ and
$\overline{A(k,l)}$ otherwise.  The expansion is uniform as $(z,w)$
varies over a compact set of strictly minimal poles with $(r,s)
\in \direc (z,w)$ and $k$ and $l_0$ not changing.
\end{theorem}

\begin{figure}
\scalebox{0.6}{\input{p1-fig2.tex}}
\vspace{10pt}
\caption{$\sing$ for Example~\ref{ex:cubic}}
\label{fig:cuberoot}
\end{figure}

\begin{example}[Cube root asymptotics] \label{ex:cubic} 

Let $F(z,w) = 1 / (3 - 3z - w + z^2)$.  The set $\sing$ is the set
$\{ w = z^2 - 3z + 3 \}$ and $g(z) = z^2 - 3z + 3$.  The point
$(1,1)$ is in $\sing$, indicating that the maximal exponential growth
rate will be zero.  Indeed, for directions above the diagonal,
Theorem~\ref{th:simple} or~\ref{th:2d} may be used at the minimal
points $\{ (z, g(z)) : 0 < z < 1 \}$, while each direction below
the diagonal corresponds to a pair of complex minimal points
fitting the hypotheses of Corollary~\ref{cor:finitely}; the result
is that the coefficients decay exponentially at a rate that is
uniform over compact subsets of directions not containing the
diagonal.

The interesting behavior is near the diagonal.  The relevant
minimal point is $(1,1)$, where $z^r w^s \equiv 1$ and the decay
is sub-exponential.  Computing $\ft''(0)$ via
equation~(\ref{eq:ft''}) below gives
$$\ft'' (z) =  -3 {z (z^2 - 4z + 3) \over (z^2 - 3z + 3)^2} \, .$$
This vanishes when $z = 1$, and computing further, we find that
$\ft$ vanishes to order exactly 3 here, with $c_3 := \ft'''(0) /
3! = i$.  Along with $\psit (0) = 1$, this then results in an
asymptotic expansion whose leading term is given by
$$a_{r,r} \sim {1 \over 2 \pi} A (3,0) i^{-1/3} (1 + e^{- i \pi /
   3}) r^{-1/3} = {\Gamma(2/3) \over 6 \sqrt{3} \pi} r^{-1/3}\, .$$
In Section~\ref{ss:open} we discuss the question of computing
asymptotics ``in the gaps'' so as to be able to conclude that
$\limsup \log a_\b{r} / \log |\b{r}| = -1/3$ or even $\limsup |\b{r}|^{1/3}
a_\b{r} = {\Gamma(2/3) \over 6 \sqrt{3} \pi}$.
 \end{example}
 
For more than two variables a result holds similar to the
two-variable result.
\begin{theorem} \label{th:higher d}
Let $F = \numer / \denom = \sum a_\b{r} \b{z}^\b{r}$ have a strictly
minimal, simple pole at $\b{z}$.  Suppose $z_d \denom_d$ does not
vanish.  If the Hessian of $\ft$ at $\b{z}$ is nonsingular, then
there is an expansion
$$a_\b{r} \sim \b{z}^{-\b{r}} \sum_{l \geq l_0} C_l r_d^{(1-d-l)/2}$$
where $l_0$ is the degree to which $\numer$ vanishes on $\sing$ near
the point $\b{z}$ . When
$\numer$ does not vanish at $\b{z}$ then $l_0 = 0$ and
$$C_0 = (2 \pi)^{(1-d)/2} \hess^{-1/2} { \numer (\b{z}) \over z_d H_d}$$
where $\hess$ is the determinant of the Hessian at $\b{z}$.
\end{theorem}

\begin{example}[Domino tilings] \label{ex:arctic} 

Random perfect tilings of planar regions by dominos have been a
subject of some interest, since the analysis by \citeasnoun{fisher;dimer}
of this
model for dimer packing uncovered an exact expression for the
partition function of the ensemble.  A generating function is
given in \citeasnoun{cohn-elkies-propp;aztec} which allowed the authors to
determine, after some cumbersome analysis, which parts of a
diamond-shaped region (a union of lattice squares approximating
the region $|x| + |y| \leq k$) were asymptotically deterministic
and which contained randomness in the limit as the edge size of
the diamond grew.

An easier analysis in the region of non-randomness is available
via Theorem~\ref{th:higher d} together with a slightly more
informative generating function than was used by
\citeasnoun{cohn-elkies-propp;aztec}.
In particular, let
$$F(x,y,z) = \sum_{t=0}^{\infty} \sum_{|r| + |s| \leq t}
   a_{r,s,t} x^r y^s z^t$$
be the generating function for the probability $a_{r,s,t}$ that
the tile covering position $(r,s)$ of a random diamond of size $t$
will be horizontal.  For brevity, we omit formal descriptions of
the diamond and its indexing.  We remark that the use of negative
indices (for each fixed $t$, the sum $\sum_{|r| + |s| \leq t}$ is
a polynomial in $x, x^{-1}, y$ and $y^{-1}$) does not require any
alterations in the theory (see \citeasnoun{cohn-pemantle;fixation} for
justification), and that the natural way to parametrize directions
is by the pair $(r/t , s/t)$ which varies over the diamond $|r/t|
+ |s/t| = 1$. From \citeasnoun{cohn-elkies-propp;aztec} or from the
generation algorithm in \citeasnoun{gessel-ionescu-propp}, one
finds
$$F(x,y,z) = {z/2 \over 1 - (x + x^{-1} + y + y^{-1}) z/2
   + z^2} \, .$$
\citeasnoun{cohn-pemantle;fixation} show that whenever $(r,s,t)$ satisfy
$$t = \sqrt{r^2 + s^2 + 2 \sqrt{r^2 + 1} \sqrt{s^2 + 1} - s} \, ,$$
then there is a smooth minimal point $(x,y,z)$ on the pole
manifold of $F$ for which $(r,s,t) \in \direc (x,y,z)$, yielding
exponential decay in the direction $(r,s,t)$.  The set of
directions so parametrized turns out to be the region between the
diamond $|r/t| + |s/t| = 1$ and the inscribed circle $(r/t)^2 +
(s/t)^2 = 1/2$.  Thus they recover the description of the region
of non-randomness as the complement of the inscribed circle.  They
also obtain descriptions of the region of fixation for related
tiling problems in which no other analysis has been carried out.
\end{example}

The extension of all of the above results to finitely minimal
points is routine.

\begin{corollary} \label{cor:finitely}
Suppose $\b{z}$ is a finitely minimal point of $\sing$ with
$\sing \cap \torus (\b{z}) = \{ \b{z}_1 , \ldots , \b{z}_n \}$.  Then
$$a_\b{r} \sim \sum_{j=1}^n E_j (\b{r})$$
where $E_j (\b{r})$ is the asymptotic expression given by the
previous theorems with $\b{z} = \b{z}_j$.  In other words, if there
are finitely many points on $\sing \cap \torus (\b{z})$, then sum the
contributions as if each were strictly minimal. 
\noproof
\end{corollary}

\begin{example}[Chebyshev polynomials] \label{ex:Tsch} 

Let $F(z,w) = 1 / (1 - 2zw + w^2)$ be the generating function for
Chebyshev polynomials of the second kind \cite{comtet;advanced}; of
course asymptotics for these are well known and easy to derive by
other means.  To use Corollary~\ref{cor:finitely}, first find the
minimal points for the direction $(r,s)$, which are $(i (\beta -
\beta^{-1})/2 , i \beta)$ for $\beta = \pm \sqrt{s-r \over s+r}$.
Computing $Q = 4 a^2 (1 - a^2)$ and summing the two contributions
then gives
$$a_{rs} \sim \sqrt{2 \over \pi} (-1)^{(s-r)/2} \left ( {2r \over
   \sqrt{s^2 - r^2}} \right )^{-r} \left ( \sqrt{s-r \over s+r}
   \right )^{-s} \sqrt{s+r \over r(s-r)}$$
when $r+s$ is even and zero otherwise, uniformly as $r/s$ varies
over compact subsets of $(0,1)$.

\end{example}

\section{Proofs of main results} \label{ss:proofs}

Half of each theorem is easy and follows directly from Cauchy's
formula
\begin{equation} \label{eq:cauchy}
a_\b{r} = \left ( {1 \over 2 \pi i} \right )^d
   \int_T \b{w}^{-\b{r} - \one} F(\b{w}) \, d\b{w}
\end{equation}
where the multi-exponent $\b{r} - \one$ means $(r_1 - 1 , \ldots ,
r_d - 1)$.  Indeed, if $\b{z}$ is a minimal point of $\sing$ then
letting $T$ approach $\torus (\b{z})$ from the inside, we see that
$|\b{z}^\b{r}| a_\b{r}$ does not increase exponentially.  If, furthermore,
 the hyperplane through $(\log |z_1| , \ldots , \log |z_d|)$ normal
to $\b{r}$ is not a support hyperplane for $\logdom$, then some $\b{x}
\in \logdom$ has $\b{x} \cdot \b{r} > (\log |z_1| , \ldots , \log |z_d|)
\cdot \b{r}$, and integrating on the torus $T(e^\b{x})$ shows that
$|\b{z}^\b{r}| a_\b{r}$ decreases exponentially.  All the work,
therefore, is in showing the converse, namely that when the
hyperplane normal to $\b{r}$ is a support hyperplane, then
$\b{z}^{-\b{r}}$ does give the right exponential order for $a_\b{r}$.
This is done by evaluating $a_\b{r}$.

Theorems~\ref{th:simple}--\ref{th:higher d} all begin with the
reduction of an iterated Cauchy integral to an oscillating
integral in one fewer dimension.

\begin{lemma} \label{lem:X}
Let $\b{z}$ be a strictly minimal simple pole of $F=\numer/\denom$.
Assume that $z_d \denom_d \neq 0$.  For a neighborhood
$\widetilde{\nbd}$ of $\zero$ in $\R^{d-1}$ define a quantity
\begin{equation}\label{eq:X}
\X := (2 \pi)^{1-d} \b{z}^{-\b{r}} \int_{\widetilde{\nbd}} \exp (- r_d \ft
  (\widehat{\bf \theta})) \psit (\widehat{\bf \theta}) \,d\widehat{\bf \theta}.
\end{equation}
Then the quantity
$$|\b{z}^\b{r}| \left | a_\b{r} - \X \right |$$
decreases exponentially as $\widetilde{\nbd}$ remains fixed and
$\b{r} \to \infty$.
\end{lemma}

\begin{proof}
For $\ee \in (0 , |z_d|)$, let $T$ be the torus $T(\b{z})$ shrunk in
the last coordinate by $\ee$, that is, the set of $\b{w}$ for which
$|w_j| = |z_j|$, $j < d$ and $|w_d| = |z_d| - \ee$.  Write
Cauchy's formula as an iterated integral
\begin{equation} \label{eq:iter}
a_\b{r} = \left ( {1 \over 2 \pi i} \right )^d
   \int_{T(\widehat{\b{z}})} \widehat{\b{w}}^{-\widehat{\b{r}} - \one}
   \left [ \int_{\C_1} w_d^{-r_d} F(\b{w}) \,
   {d w_d \over w_d} \right ] \, d \widehat{\b{w}}  \, .
\end{equation}
Here $\C_1$ is the circle of radius $|z_d| - \ee$.  Let $K
\subseteq T(\widehat{\b{z}})$ be a compact set not containing
$\widehat{\b{z}}$.  For each fixed $\widehat{\b{w}} \in K$, the function
$F(\widehat{\b{w}} , \cdot)$ has radius of convergence greater than
$|z_d|$.   Hence the inner integral in equation~(\ref{eq:iter})
is $O(|z_d|+ \delta)^{-r_d}$ for some $\delta > 0$.  By continuity
of the radius of convergence,we may integrate over $K$ to see that
$$|\b{z}^\b{r}| \int_{K \times \C_1} \b{w}^{-\b{r} - \one} F(\b{w}) \, d\b{w}$$
decreases exponentially.  Thus if $\nbd$ is any neighborhood of
$\widehat{\b{z}}$ in $T(\widehat {\b{z}})$, the quantity
$$|\b{z}^\b{r}| \left | a_\b{r} - \left ( {1 \over 2 \pi i} \right )^d
   \int_\nbd \widehat{\b{w}}^{-\widehat{\b{r}} - \one}
   \left [ \int_{\C_1} {F(\b{w}) \over w_d^{r_d + 1}} \, dw_d
   \right ]  d \widehat{\b{w}} \right |$$
decreases exponentially.  Thus we have reduced the problem to an
integral over a neighborhood of $\widehat{\b{z}}$.

Near $\b{z}$ there is a parametrization $w_d = g(\widehat{\b{w}})$ of
$\sing$.  Let $\C_2$ be the circle of radius $|z_d| + \ee$.  Then
when $\nbd$ is sufficiently small compared to $\ee$, the image of
$\nbd$ under $g$ is disjoint from $\C_2$.  Fix such a neighborhood.
For any $\widehat{\b{w}} \in \nbd$, the function $F(\widehat{\b{w}} ,
\cdot)$ has a single simple pole in the annulus bounded
by $\C_1$ and $\C_2$, occurring at $g(\widehat{\b{w}})$.  The
residue in the last variable of $F$ at $g(\widehat{\b{w}})$ is equal
to
\begin{equation} \label{eq:defres}
R(\widehat{\b{w}}) := - \psi (\widehat{\b{w}}) g(\widehat{\b{w}})^{-r_d}
\end{equation}
where $\psi$ is defined in~(\ref{eq:psi}). Therefore, for each
fixed $\widehat{\b{w}} \in \nbd$,
$$\int_{\C_1} {F(\b{w}) \over w_d^{r_d + 1}} \, d w_d
   = \int_{\C_2} {F(\b{w}) \over w_d^{r_d + 1}} \, d w_d
   - 2 \pi i R(\widehat{\b{w}}) .$$
But $|\b{z}^\b{r} \int_{\C_2} F(\b{w}) dw_d / \b{w}^{\b{r} + \one}|$ is
bounded by a constant multiple of $(1 + \ee / |z_d|)^{-r_d}$ (the
constant depending on the maximum of $F$ on $\C_2$) and hence
$|\b{z}^\b{r}| |a_\b{r} - X|$ is exponentially decreasing, where
\begin{eqnarray} \label{eq:res}
X & = & (2 \pi i)^{1-d} \int_\nbd (\widehat{\b{w}})^{-\widehat{\b{r}}
   - 1} g(\widehat{\b{w}})^{-r_d} \psi (\widehat{\b{w}}) \,
   d \widehat{\b{w}} \\[1ex]
& = & (2 \pi i)^{1-d} \b{z}^{-\b{r}} \int_\nbd {\widehat{\b{w}}^{-\widehat{\b{r}}}
   \over \widehat{\b{z}}^{-\widehat{\b{r}}}} {d \widehat{\b{w}} \over
   \prod_{j=1}^{d-1} w_j} \left ( {g(\widehat{\b{w}}) \over g(z_d)}
   \right )^{-r_d} \psi (\widehat{\b{w}}) \nonumber
\end{eqnarray}
Changing variables to
$w_j = z_j e^{i \theta_j}$ and $dw_j = i w_j d\theta_j$ turns the
quantity $X$ into
$$(2 \pi)^{1-d} \b{z}^{-\b{r}} \int_{\widetilde{\nbd}} \prod_{j=1}^{d-1}
   e^{-i r_j \theta_j} \psit (\widehat{\bf \theta}) \left ( {
   g(\widehat{\b{w}}) \over g(\widehat{\b{z}})} \right )^{-r_d}
   \, d\widehat{\bf \theta} $$
and plugging in the definitions of $f$ and $\ft$ at~(\ref{eq:f})
above yields
$$(2 \pi)^{1-d} \b{z}^{-\b{r}} \int_{\widetilde{\nbd}} \exp (- r_d
   \ft (\widehat{\bf \theta})) \psit (\widehat{\bf \theta}) \,
   d\widehat{\bf \theta}$$
which is none other than $\Xi$.   \end{proof}

\begin{unremark} It is possible to compute from Cauchy's
integral formula in a more coordinate-free way as follows.  There
is a unique holomorphic $(d-1)$-form $\omega_F$ on $\sing$ for which
$\omega \wedge d\denom = \numer \, dz_1 \wedge \cdots \wedge dz_d$.
Let $\Omega$ be a ($d+1$)-manifold that is a homotopy from a small
torus to a torus at infinity.  Then $M := \Omega \cap \sing$ is a
$(d-1)$-manifold and $a_\b{r} = (2 \pi i)^{-d} \int_M \b{w}^{-\b{r} - \one} dF$
in the sense of currents, which is none other than $\int_M
\b{w}^{\b{r} - \one} \omega_F$.  See \citeasnoun{kenyon-pemantle} for a
more thorough discussion of the foregoing.  The manifold $M$ is
any member of a certain homology class in $\sing$ with the coordinate
axes removed, and choosing $M$ to pass through the stationary
phase point for the integrand replicates the selection of $\b{z}$
with $\b{r} \in \direc(\b{z})$.  Although more canonical, the
coordinate-free method is less suitable for explicit computation,
so we do not pursue it further here.  Suffice it to point out that
the conclusion of Theorem~\ref{th:higher d} may of course be
written in terms more evidently symmetric, as was done in
Theorem~\ref{th:simple}.
\end{unremark}

Equation~(\ref{eq:X}) is easily recognized as the standard form
for an oscillating integral.  The only unusual feature is that the
phase is neither real nor purely imaginary.  This presents no
difficulties, but it does necessitate the statement of a result in
Section~\ref{ss:oscillate} that is a little different from the
usual results on purely oscillating integrals, found in, for example,
\citeasnoun{stein;oscillatory-integrals} or
\citeasnoun{bleistein-handelsman;asymptotic-expansions-integrals}.  We first
establish that $\widehat{\bf \theta} = \zero$ is a stationary
phase point for the function $\ft$ when $\b{r} \in \direc (\b{z})$.
\begin{lemma} \label{lem:stat}
The quantity $\ft (\zero)$ always vanishes.  If $\b{r} \in \direc
(\b{z})$ then $\grad \ft (\zero) = \zero$ and the real part of $\ft$
has a strict minimum at $\zero$.
\end{lemma}

\begin{proof} The first statement is immediate.  To prove
the second, let $j \leq d-1$ and see from the definition of $f$
that
$$r_d f_j (\widehat{\b{z}}) = {r_d g_j (\widehat{\b{z}}) \over
   g(\widehat{\b{z}})} + {r_j \over z_j}.$$
By definition of $\direc$, the ratio $r_j / (z_j \denom_j)$ is
some constant $c$ independent of $j$, hence
$$c^{-1} r_d f(\b{z}) = g_j (\b{z}) \denom_d (\b{z}) + \denom_j (\b{z}) .$$
The right hand side of this is the derivative of $\denom (w_1 ,
\ldots , w_{d-1} , g(\widehat{\b{w}}))$ with respect to $w_j$ at
$\widehat{\b{z}}$.  By definition of $g$ this vanishes, and hence
$f_j (\widehat{\b{z}}) = 0$.  But $\ft_j (\zero) = i z_j f_j (\b{z})$, so
the gradient of $\ft$ must vanish at $\zero$.  Finally, observe
that $\re \{ \ft (\widehat{\bf \theta}) \} = -\log |\gt
(\widehat{\bf \theta}) / z_d|$.  By strict minimality of $\b{z}$,
the modulus of $g(\widehat{\b{w}}) = \gt (
\widehat{\bf \theta})$ is greater than $|z_d|$ for any
$\widehat{\b{w}} \in \torus (\widehat{\b{z}})$.   \end{proof}

We now prove Theorems~\ref{th:simple},~\ref{th:2d}
and~\ref{th:higher d} in reverse order.  We see from
Lemma~\ref{lem:X} that proving any of these theorems amounts to
evaluating the quantity $\X$ in equation~(\ref{eq:X}).  From
Lemma~\ref{lem:stat} we see that $\zero$ is a stationary point for
the function $\ft$ as long as $\b{r} \in \direc (\b{z})$.  The
function $\ft$ is in general complex valued, but we will see in
Theorem~\ref{th:nondeg} that it may be treated as if it were real
valued, given the strict minimality of the zero guaranteed by
Lemma~\ref{lem:stat} and the nonsingularity hypothesis. In
particular the leading term of the integral in~(\ref{eq:X}) is
$(2 \pi)^{(d-1)/2} \psit (\zero) r_d^{(1-d)/2}$ divided by the
product of the square roots of the eigenvalues of the Hessian.
Once we have identified $\psit (\zero) = \psi (\zero)$ as $\numer
(\zero) / (z_d \denom_d)$, the theorem follows directly from
Theorem~\ref{th:nondeg}.

Theorem~\ref{th:2d} follows from the more explicit asymptotic
development given in Corollary~\ref{cor:full 2d}. Finally, to
prove Theorem~\ref{th:simple}, it remains to compute the quantity
$\ft'' (0)$ in terms of the partial derivatives of $\denom$.
First we compute the derivatives of $g$.

\begin{lemma} \label{lem:der g}
In a neighborhood of $(z,w)$, $\psi$ and the derivatives of $g$
are as follows.
\begin{eqnarray}
g'(z) & = & -{\denom_z \over \denom_w} \label{eq:g'} \\[2ex]
g''(z) & = & - {1 \over \denom_w}  \left [ \denom_{zz} - 2
   {\denom_z \over \denom_w} \denom_{zw} + {\denom_z^2 \over
   \denom_w^2} \denom_{ww} \right ] . \label{eq:g''} \\[2ex]
\psi (z) & = & {\numer (z,w) \over - w \denom_w(z,w)} . \nonumber
\end{eqnarray}
\end{lemma}

\begin{proof} Differentiate the equation $\denom(z,g(z)) = 0$
to get $\denom_z + g'(z) \denom_w = 0$ which is the same
as~(\ref{eq:g'}).  Differentiate again to get
$$\denom_{zz} + 2 g' \denom_{zw} + g'' \denom_w + (g')^2 \denom_{ww} = 0$$
and use~(\ref{eq:g'}) to eliminate $g'$, giving~(\ref{eq:g''}).
The formula for $\psi$ follows from the definitions of $\psi$ and
of the partial derivative.
\end{proof}

\noindent{\sc Proof of Theorem}~\ref{th:simple}~{\sc via direct
computation:}  We know from Lemma~\ref{lem:stat} that $\ft$
vanishes to order at least two at 0.  To compute $\ft''(0)$,
observe first that $\ft'' - \log \gt$ is linear in $\theta$, so
$\ft'' = (\log \gt)''$.  When $Z = z e^{i \theta}$, we
have $(d/d\theta) = iZ (d/dZ)$, so
$$\ft'' = i Z {d \over dZ} \left ( iZ {d \log g \over dZ} \right )
   = - Z {d \over dZ} \left ( {Z g' \over g} \right ) \, .$$
Expanding this yields
\begin{equation} \label{eq:ft''}
\ft'' = - Z {g' + Z g'' \over g} + {Z^2 (g')^2 \over g^2} \, .
\end{equation}

By our assumption, $G$ does not vanish at $(z,w)$, so as long as
$\ft'' (0) \neq 0$, we may use Theorem~\ref{th:2d} to
conclude that the leading term asymptotic for $a_{r,s}$ is the
$k=2, l=0$ term of~(\ref{eq:full 2d}).  The term $b_0$ there is
equal to
$$\psit(0) \eta' (0) = \psi (z) \sqrt{2 / \ft''(0)} =  {\numer (z,w)
   \over - w \denom_w (z,w)} \sqrt{2 \over \ft'' (0)} .$$
Thus from Theorem~\ref{th:2d},
$$a_{r,s} \sim {A (2,0) \over 2 \pi} z^{-r} w^{-s} {\numer (z,w)
   \over w \denom_w (z,w)} \sqrt{2 \over s \ft'' (0)} .$$
Now evaluate this using the value $A (2,0) = \sqrt{\pi}$ and
equation~(\ref{eq:ft''}) along with~(\ref{eq:g'})
and~(\ref{eq:g''}) to obtain
$$a_{r,s} \sim {1 \over \sqrt{2 \pi}} z^{-r} w^{-s} {\numer (z,w)
   \over w \denom_w (z,w)} \sqrt{(-w \denom_w (z,w))^3 \over s Q}$$
where
$$Q = (- w \denom_w (z,w))^3 \ft''(0) = (- w \denom_w (z,w))^3
   z {-g'(z) - z g''(z) \over g(z)} + {z^2 (g'(z))^2 \over (g(z))^2} \,
   .$$
With the help of Lemma~\ref{lem:der g} we see (using $g(z) = w$)
that
$$Q = (- w \denom_w)^3 \left [ -z {\denom_z \over -w
   \denom_w} - z^2 {1 \over -w \denom_w}
   \left (\denom_{zz} - 2 {\denom_z \over \denom_w} \denom_{zw}
   + {\denom_z^2 \over \denom_w^2} \denom_{ww} \right )
   + {z^2 \denom_z^2 \over w^2 \denom_w^2} \right ] \, ,$$
evaluated at $(z,w)$, which simplifies to the expression in
Theorem~\ref{th:simple}.  We see also that the nonvanishing
hypotheses on $Q$ is enough to guarantee $\ft''(0) \neq 0$, which
finishes the proof of Theorem~\ref{th:simple}.
\noproof

\section{Some oscillating integrals} \label{ss:oscillate}

The oscillating integrals we require are integrals over a
neighborhood of zero in $\R^d$ of the complex-valued integrand:
$$\int_\nbd \exp (- \lambda f (\b{x})) \psi (\b{x}) \, d\b{x}$$
where $f(\zero) = 0, \grad f (\zero) = \zero$ and $\re \{ f
\} \geq 0$.  They are not difficult to compute, but since the
standard references assume $f$ is either real or purely imaginary,
we sketch the development of these results.  We mostly follow the
exposition of \citeasnoun{stein;oscillatory-integrals}, adapting
it to complex-valued phase functions and simplifying it to take
advantage of the decay of the magnitude of the integrand in this
case.

We begin with one-dimensional results.  Let $\cio$ denote the
class of smooth functions with compact support.  The following
proposition is a well known consequence of Watson's Lemma (see,
for example, \citeasnoun[Ch.~2, Theorem~1]{wong;asymptotic-integrals}).
\begin{proposition} \label{pr:laplace}
Let $\psi \in \cio (\R)$ and denote $b_j = \psi^{(j)} (0) / j!$.
Then as $\lambda \to \infty$, there is an asymptotic development
$$\int_0^\infty \exp (- \lambda x^k) \psi (x) \, dx \sim
   \sum_{l=0}^\infty A_+ (k,l) b_l \lambda^{-(l+1)/k},$$
where, as in~(\ref{eq:plus}),
$$A_+ (k,l) := k^{-1} \Gamma \left ( {l+1 \over k} \right ) .$$
\noproof
\end{proposition}

We extend this to more general one-sided integrals by a complex
change of variables. Given any analytic, complex-valued function
$f$ on an interval $[0,B]$, suppose that $f(0) = 0$, that
$f' \neq 0$ on $(0,B]$, and let $k \geq 1$  be the minimal so that
$f^{(k)} (0) \neq 0$.  Let $\psi \in \cio$ vanishing to order $l
\geq 0$ at 0.  Denote $c_j = f^{(j)} (0) / j!$ and $b_j =
\psi^{(j)} (0) / j!$.  The real part of $c_k$ is necessarily
nonnegative.  Define a function $y$ on $[0,B]$ by
$$y(x) = f(x)^{1/k} = c_k^{1/k} x \left ( 1 + {f(x) - c_k x^k
   \over c_k x^k} \right )^{1/k} \, ,$$
where the argument of $c_k^{1/k}$ is between $-\pi/(2k)$ and $\pi
/ (2k)$.  The quantity $f(x) - c_k x^k$ is $O(x^{k+1})$ near zero,
so $y$ is analytic near 0, and, in particular, is a diffeomorphism
between $[0,B]$ and a contour $\gamma$ from 0 to some $B^*$.  Let
$F$ invert $y$. The derivatives of
$F$ at 0 are easy to compute formally and the first $j+1$ starting
from the $k^{th}$ depend only on the first $j$ coefficients of $f$
starting at $c_k$. Define
\begin{eqnarray}
\psi^* & = & (\psi \circ F) \cdot F' \; ; \nonumber \\
b^*_j & = & \psi^*\,^{(j)} (0) / j! \; . \label{eq:bs}
\end{eqnarray}
\begin{theorem} \label{th:one-sided}
Let $f$ be analytic (complex-valued) on an interval
$[0,B]$. Assume that $f(0) = 0$, that $f' \neq 0$ on $(0,B]$,
and $\re \{ f \}$ has a strict minimum at 0.  Let $k \geq 2$ be
minimal such that $f^{(k)} (0) \neq 0$ and $m$ be minimal so that
the real part of $f^{(m)} (0)$ does not vanish. Let $\psi \in
\cio$, let $l$ be minimal such that $\psi^{(l)} (0) \neq 0$,
and denote $c_j := f^{(j)} (0) / j!$, $b_j := \psi^{(j)} (0) /
j!$. Define $b^*_j$ as in~(\ref{eq:bs}).  Then there is an
asymptotic development
\begin{equation} \label{eq:series}
\int_0^B \exp (- \lambda f(x)) \psi (x) \, dx \sim
   \sum_{j=l}^\infty A_+ (k,j) b^*_j \lambda^{-(j+1)/k} .
\end{equation}
The constant in the $O(\lambda^{-(N+1)/k})$ term depends continuously (only) on the derivatives of $f$ and $\psi$ up to $(N+1) m / k - 1$.
\end{theorem}

\begin{proof} Changing variables to $y = f(x)^{1/k}$,
the integral becomes
$$\int_\gamma \exp (- \lambda y^k) \psi^* (y) \, dy ;$$
the curve $\gamma$ is the image of $[0,B]$ under $y$, so
$\gamma' (0) = c_k^{1/k}$ and $\gamma$ remains in the
right half plane, strictly except at 0.  For $0 < N < M$ write
$\psi^*$ as $P_M + y^{M+1} R_M$, where $P_M$ is a polynomial of
degree $M$ and $R_M$ is bounded; this can be done since $\psi^*$
may be approximated by a degree $M$ polynomial to within
$O(y^{M+1})$ at 0.

First, evaluate
$$\int_\gamma \exp (- \lambda y^k) P_M (y) \, dy$$
by moving the contour.  Replace $\gamma$ by two line segments, the
first of which goes along the positive real axis to some distance
$\ee$ and the second of which is strictly in the right half plane
(we assumed $\re \{ f \} > 0$ except at 0). The integral along the
second segment is exponentially small since the integrand is.
Hence the combined contribution is the series~(\ref{eq:series})
out to the $j=M$ term.

Next, bound
$$\left | \int_\gamma \exp (- \lambda y^k) y^{M+1} R_M (y) \, dy
   \right | \, .$$
With $C$ representing different constants in different lines, we
now observe that on $\gamma$ we have $\re \{ - y^k \} < - C
|y|^m$. Thus, parametrizing $\gamma$ by arc-length, an upper bound
is given by
$$\int_0^\infty \exp (- \lambda C t^m) t^{M+1} C |R_M (\gamma
   (t))| \, dt .$$
This is easily seen to be bounded above by
$C \lambda^{-(M+2)/m}$ where $C$ depends on the first $M$
derivatives of $f$ and $\psi$.  Choosing $M \geq m(N+1)k - 1$ we
have a remainder term that is $O(\lambda^{-(N+1)/k})$, proving the
theorem. 
\end{proof}

The value of a two-sided integral follows as a corollary.
\begin{corollary} \label{cor:full 2d}
Assume the hypotheses of Theorem~\ref{th:one-sided}, with $f$ now
defined on an interval $[-B,B]$.  Then there is an asymptotic
development
\begin{equation}\label{eq:two-sided}
\int_{-B}^B \exp (- \lambda f(x)) \psi (x) \, dx \sim
   \sum_{j=l}^\infty A (k,j) b^*_j \lambda^{-(j+1)/k} .
\end{equation}
with $A (k,j)$ given by~(\ref{eq:k odd}) and~(\ref{eq:k even}).
The bounds on the remainder terms each depend continuously on finitely many
derivatives of $f$ and $\psi$ on $[-B,B]$.
\end{corollary}

\begin{proof} The two-sided integral is the sum of two
one-sided integrals on intervals $[0,B]$ and $[-B,0]$.  The
integral over $[-B,0]$ may be written as an integral over $[0,B]$
of the function $\exp (- \lambda f(-x)) \psi (-x) \, dx$. With
$b_l^*$ still denoting the coefficients resulting form the application
of Theorem~\ref{th:one-sided} to the first integral, let
$\check{b_l^*}$ denote the coefficients when
Theorem~\ref{th:one-sided} is applied to the second integral.  In
order to add the two integrals, we write $\check{b_l^*}$ in terms
of $b_l^*$ by means of the following routine computation.

Let $c_k := c^k e^{i \alpha}$ with $c > 0$ and $|\alpha| \leq \pi
/ 2$ and define the analytic quantity $R$ so that
$$y(x) = \left [ c_k x^k (1 + R(x))^k \right ]^{1/k} =
   c e^{i \alpha / k} x (1 + R(x)) .$$
If $k$ is odd, then then the hypothesis $\re \{ f
\} \geq 0$ implies that $c_k$ is purely imaginary.  We have
$$\check{y} (x) = \left [ - c_k x^k (1 + R(-x))^k \right ]^{1/k} =
   c e^{- i \alpha / k} x (1 + R(-x)) = - y(-x) e^{- 2 i \alpha / k}
   .$$
Writing $\finv$ for the inverse function to $y$ and $\check{\eta}$
for the inverse function to $\check{y}$ we then have
$$\check{\finv} (x) = - \finv (-e^{2 i \alpha / k} x) \, .$$
Hence, letting $\C_l [ \cdot ]$ denote the coefficient of
$y^l$,
\begin{eqnarray*}
\check{b_l^*} & = & \C_l \left [ \psi ( - \check{\finv} (x))
   \cdot \check{\finv}' (x)\right ] \\[1ex]
& = & \C_l \left [ \psi ( \finv (-e^{2 i \alpha / k} x))
   \cdot e^{2 i \alpha / k} \cdot \finv' (-e^{2 i \alpha / k} x) \right ]
\end{eqnarray*}
and thus
$$\check{b_l^*} = (-1)^l e^{2 i \alpha (l+1) / k} b_l^* \, .$$
When $k$ is even, the computation is similar but easier, resulting
in
$$\check{b_l^*} = (-1)^l b_l^* \, .$$
Now observe that if $k$ is odd, hence $c_k$ is purely imaginary,
then $e^{2 i \alpha (l+1) / k} = e^{\pm i \pi (l+1) / k}$
according to the sign of the argument of $c_k$.  Setting $A (k,l)
= (1 + (-1)^l) A_+ (k,l)$ if $k$ is even and
$(1 + e^{\sign\Arg (c_k) i \pi (l+1) / k} A_+ (k,l)$
if $k$ is odd, we recover the definition in~(\ref{eq:k odd})
and~(\ref{eq:k even}) and prove the Corollary.
\end{proof}
\begin{theorem} \label{th:nondeg}
Let $f$ be a smooth complex-valued function on a neighborhood of
$\zero$ in $\R^d$ such that $\re \{ f \} \geq 0$ with equality only at $\zero$.  Suppose further that $\grad f(\b{0})=0$, and that the Hessian (matrix of second partials) of $f$ has eigenvalues with positive real parts. Let $\hess$ denote the Hessian determinant at 0. Then for $\psi \in \cio$,  there is an asymptotic expansion
$$\int \exp (- \lambda f (\b{x})) \psi (\b{x}) \, d\b{x} \sim
   \sum_{j \geq l} C_j \lambda^{-(l+d)/2}$$
where $l$ is the degree of vanishing of $\psi$ at $\zero$. If $l =
0$ then $C_0 = \psi (\zero) (2 \pi)^{d/2} \hess^{-1/2}$. The
choice of square root is determined by $\hess^{-1/2} =
\prod_{j=1}^d \mu_j^{-1/2}$ where $\mu_j$ are the eigenvalues of
the Hessian and the principal square root is taken in each case.
\end{theorem}

\begin{proof} Let $Q = \sum_{i,j=1}^d q_{i,j} z_i z_j$
be the quadratic form determined by the Hessian at the origin.
Denote the eigenvalues of $Q$ by $\{ \mu_j : 1 \leq j \leq d \}$
and note that each $\mu_j$ has nonnegative real part.

Step 1: change coordinates to make $f$ exactly equal to the
quadratic form $Q$. Indeed since $f(\b{x}) = Q(\b{x})/2 + O(|\b{x}|^3)$,
and the Hessian is nondegenerate, there is a locally smooth change
of variables $\{ x_j (\b{z}) : 1 \leq j \leq d \}$ such that $f(\b{z})
= Q(\b{x} (\b{z}))/2$ and the Jacobian at the origin is 1.

Step 2: normalize by $\hess^{1/2}$.  For any quadratic form $Q$
there is a linear change of variables $\b{y} (\b{x})$ such that
$Q(\b{x}) =
\sum_{j=1}^d y_j^2$.  The change of variables matrix $P$
satisfies $P P^T = M(Q)$, the symmetric matrix representing $Q$.
Changing variables to $\b{y}$ introduces an integrating factor of
$\det P$ which is a square root of $\hess$ since $M(Q)$ is just the
Hessian.  Let $\nbd'$ be the region of integration over which
$\b{y}$ varies when $\b{z}$ varies over an appropriately small
neighborhood of $\zero$.

Step 3: Expand $\psit$ into monomials.  The function $\psi$ has
now become $\psit$, where $\psit(\zero) = \hess^{-1/2} \psi
(\zero)$ and the sign of the square root will be chosen later.  We
may expand $\psit$ into monomials, using the same argument as in
the proof of Theorem~\ref{th:one-sided} to show the remainder term
can be made
$O(|\b{y}|^N)$ for any $N$.  It remains to evaluate the integral over the
region of integration, $\nbd'$ of
$$\int_{\nbd'} \exp (- \lambda \sum_{j=1}^d y_j^2) \psit (\b{y}) \, d\b{y}$$
when $\psit$ is a monomial.

Step 4: move the region of integration to the real $d$-space. Let
$\nbd''$ be the projection of $\nbd'$ onto $\R^d$ by setting the
imaginary part to zero.  We claim that changing the region of
integration from $\nbd'$ to $\nbd''$ alters the integral by an
amount rapidly decreasing in $\lambda$.  To show this, let
$\Omega$ be the region $\{ \re \{ \b{x} \} + i t \im \{ \b{x} \} :
\b{x} \in \nbd' , t \in [0,1] \}$.  The boundary of $\Omega$
(as a manifold) is composed of $\nbd', \nbd''$
(with opposite signs) together with $S := \{ \re \{ \b{x} \} + i t
\im \{ \b{x} \} : \b{x} \in \partial \nbd' , t \in [0,1] \}$.  For any
$d$-form $\omega$, $\int_\Omega d\omega = \int_{\partial \Omega}
\omega$.  When $\omega = \exp (- \lambda \sum_{j=1}^d \mu_j y_j^2)
\b{y}^\b{r} dy_1 \wedge \cdots \wedge dy_d$ is a holomorphic $d$-form,
we see that $d\omega$ vanishes (being the sum of $\partial /
\partial \overline{z_j}$ terms) so that
$$\int_{\nbd'} \omega = \int_{\nbd''} \omega + \int_S \omega .$$
We know that $\re \{ \sum_j \mu_j y_j^2 \}$ is bounded away from 0
on $\partial \nbd'$, and its minimal value on $S$ lies on
$\partial \nbd'$, hence the integral over $S$ decays
exponentially.

Step 5: evaluate the integral.  Factoring $\int_{\nbd''} \b{y}^\b{r}
\exp (- \lambda \sum_{j=1}^d y_j^2)$ into one-dimensional
integrals and plugging into Proposition~\ref{pr:laplace} yields an
asymptotic expansion whose leading term (when $l=0$) is equal to
$(2 \pi)^{d/2} \psi (\zero) \hess^{-1/2}$.  When $f (\b{z})$ is the function $\sum_{j=1}^d z_j^2$, then the positive square root is taken.  The
choice of square root must be continuous in the analytic topology
on functions having nondegenerate Hessians and having eigenvalues
with positive real parts, and the only such choice is the product
of the principal square roots of the eigenvalues of the Hessian.
\end{proof}

\section{Classification of cases} \label{ss:classify}

For purposes of classification some natural questions are:
\begin{enumerate}
\romenumi
\item what are all possible local geometries of minimal points of $\sing$?
\item which of these can be handled by variants of the methods in this paper?
\item are these sufficient to yield a good approximation to $a_\b{r}$ no
matter what the direction, $\b{r} / |\b{r}|$, and no matter which
generating function in the class, say, of functions meromorphic in
a neighborhood of their domain of convergence?
\end{enumerate}
To make the last question more concrete, consider the simplest
possible example, namely binomial coefficients, where $F =
1/(1-z-w)$ and $\sing$ is a complex line.  There are no singular
points here, but how do we know that as $(z,w)$ varies over
minimal points of $\sing$, the direction $\direc (z,w)$ will cover
all of $\RP^1$?

This question will be answered by Theorem~\ref{th:classify}, but
first we need to add some detail to the geometric discussion begun
in Section~\ref{ss:geom}.  It will be evident that quite a few
cases need to be considered, some of which require new tools and
some of which require only minor modifications.  Accordingly, the
results will appear in several papers, currently under
preparation.  In other words, a discussion of taxonomy will
necessarily refer to results not yet published, and we will
indicate to the best of our knowledge which ones are expected to
be routine.

Given a point $\b{z} \in \sing$, we extend the definition of
$\direc (\b{z})$ to mean the set of limits of $\direc (\b{y})$ as $\b{y}
\to \b{z}$ along smooth points.  When $\b{z}$ is minimal, this is just
the set of normals to support hyperplanes of $\logdom$ at the
point $(\log |z_1| , \ldots , \log |z_d|)$, so this is consistent
with the old definition.  As we will see shortly, $\direc (\b{z})$
may be a $(d-1)$-dimensional subset of $\RP^{d-1}$ when $\b{z}$ is a
critical point of $\sing$.

When $\denom$ has a repeated factor, the residue computation in
equation~(\ref{eq:defres}) must be replaced by one involving the
derivative.  The remainder of the computation proceeds without a
hitch as before.  Details are given in \citeasnoun{pemantle-wilson;toral}.
For the remainder of the taxonomy, we assume $\denom$ to
be square-free.  Toral smooth points may be handled by methods
exactly the same as strictly minimal points.  The inner integrand
in~(\ref{eq:iter}) will in this case have its maximal modulus on a
set of dimension larger than zero.  A modification of the
necessary oscillating integral computation that works in this case
is also given in \citeasnoun{pemantle-wilson;toral}.

If $\b{z} \in \sing$ is not smooth, all the first partials vanish.  The
expansion of $\denom (\b{x})$ near $\b{z}$ is then a sum of terms of
degrees 2 and higher.  We call $\b{z}$ a {\em homogeneous point} of
degree $k$ if this expansion contains terms
$(x_j - z_j)^k$ for each $j = 1 ,\ldots , d$, and contains no
terms of total degree less than $k$.

\begin{lemma} \label{lem:homogeneous}
If $\b{z}$ is a locally minimal point of $\sing$ with nonzero
coordinates, and $F$ is meromorphic in a neighborhood of $\b{z}$
then $\b{z}$ is homogeneous.
\end{lemma}

\begin{proof} Passing to $F(z_1 x_1 , \ldots , z_d x_d)$ if
necessary, we may assume $\b{z} = \one$.  Setting $x_j = 1$ for all
but one index $j$, we cannot obtain the zero function (by
minimality), and so some term in the expansion around $\one$ is a
pure power of $(x_j - 1)$, and we denote the minimal degree such
term by $c_j (x_j - 1)^{k_j}$.  If $\b{z}$ is not a homogeneous
point, then there is some $j$ for which some monomial has total
degree lower than $k_j$.  Assume without loss of generality that
$j=d$.  The function $F(x , x  , \ldots , x , y)$ then has a
minimal degree pure $y-1$ term $c_0 (y-1)^k$,
$k := k_d$, and some term $c' (x-1)^a (y-1)^b$ with $a+b < k$.
In other words, the Newton Polygon of $F(x, \ldots , x , y)$
around $(1,1)$ has a support line passing through $(0,k)$ with
slope $-p/q$ in lowest terms, and $p > q$.  It is well known that
we may describe the solutions $y(x)$ of the equation
$$F(1+x , \ldots , 1+x , 1 + y) = 0$$
as follows.  Write
$$H := (y-1)^k (c_0 + c_1 (y-1)^{-p} (x-1)^q + c_2(y-1)^{-2p} (x-1)^{2q}
  + \cdots + c_s (y-1)^{-sp} (x-1)^{sq})$$
for the polynomial collecting all the terms on this support line.
Then for each $q^{th}$ root of unity, $\omega$, and each root
$\lambda$ of $\sum c_{s-j} \lambda^j = 0$, there is a solution
$y = \lambda^{1/p} x^{q/p} (\omega + o(1))$ as $x \to 0$.
A proof may be found in \citeasnoun{brieskorn-knorrer}.

Varying $x$ over the set $|\pi - \arg (x)| \leq \pi / 4$, we see
that the solutions $y(x)$ must sometimes be in this set as well.
For those $x$, the points $(1+x , \ldots , 1+x , 1+y)$ will be in
$\sing \cap \disk (\one) \setminus \torus (\one)$, violating
minimality of $\one$.  By contradiction, we have shown that no
monomial in the expansion around $\one$ has lower total degree
than any pure power term, hence $\one$ is minimal.    \end{proof}

Continuing the taxonomy, suppose that $\b{z}$ is a homogeneous point
of $\sing$ of degree $k \geq 2$.  We say that $\b{z}$ is a {\em
multiple point} if $\sing$ is locally the union of $k$ analytic
surfaces.  Algebraically, this means that the leading (order $k$)
terms in the expansion of $\denom$ near $\b{z}$ factors into linear
pieces.  If the homogeneous point $\b{z}$ is not a multiple point,
we say it is a {\em cone point}.  When $d=2$ there are no cone
points, since any homogeneous polynomial in 2 variables factors
completely over $\CC$.

Our understanding of cone points is not yet complete, but an
analysis involving cone points is underway in \citeasnoun{cohn-pemantle;fixation}.
For multiple points, most of the story is given in
\citeasnoun{pemantle-wilson;multiple}. In particular, the following theorem
is proved there.
\begin{theorem}[\citeasnoun{pemantle-wilson;multiple}]
\label{th:PW2000a}
Let $\b{z}$ be an isolated, minimal, multiple point of $\sing$ with
multiplicity $k$.  Let $S \subseteq \RP^{d-1}$ be the set of
outward normals to support hyperplanes to $\logdom$ at the point
$(\log |z_1| , \ldots , \log |z_d|)$.  Then there is an integer $p \geq 0$
and a polynomial function $\phi : S \to \R$ such that the
asymptotic expansion
\begin{equation}\label{eq:PW2000a}
a_\b{r} \sim \b{z}^{-\b{r}} \phi (\b{r}) \sum_j C_j (r_d)^{k-p/2-j/2}
\end{equation}
holds uniformly as $\b{r}$ varies over compact subsets of the
interior of $S$.
\noproof
\end{theorem}

The extension to toral multiple points is given in
\citeasnoun{pemantle-wilson;toral}. If a multiple point is not
isolated or toral, then the degree of multiplicity, $k$, must be
less than the dimension,
$d$.  This cannot happen of course when $d=2$, but does happen
when $d \geq 3$.  The method for handling this case, toral or
otherwise, is given in \citeasnoun{pemantle-wilson;toral}.  That paper
will also contain some subcases of the isolated multiple point
case, namely when the sheets of $\sing$ intersect non-transversely.

Having more or less completed the taxonomy, we now discuss when we
can guarantee that our methods yield asymptotics in all
directions.
\begin{theorem} \label{th:classify}
Let $F = \numer / \denom = \sum a_{r,s} z^r w^s$ be the quotient
of analytic functions $\numer , \denom : \CC^2 \to \CC$.  Suppose
that the coefficients $a_{r,s}$ are all nonnegative, and that
$F(z,0)$ and $F(0,w)$ are not entire.  Then for every direction
$\alpha \in \RP^1$ there is a minimal $\b{z} \in \sing$ with
$\alpha \in \direc (\b{z})$.
\end{theorem}

\begin{proof} Let $(x,y)$ be any point on the boundary of
$\logdom$.  For $u < e^x$ and $v < e^y$ the power series for $F$
is convergent at $(u,v)$.  As $u \uparrow e^x$ and $v \uparrow
e^y$ therefore, $F(u,v)$ is finite and increasing.  On the other
hand, the power series for $F$ is not absolutely convergent on
$\torus (e^x , e^y)$, since we know $F$ to have some singularity
on this torus.  Hence $F(u,v) \uparrow \infty$ as $(u,v) \uparrow
(e^x,e^y)$.  Since $F$ is meromorphic, it must have a pole at
$(e^x , e^y)$, hence $(e^x , e^y) \in \sing$ and is a minimal point of $\sing$.
As $(x,y)$ varies over the boundary of $\logdom$, we let $\gamma
\subseteq \sing$ denote the curve traced out by this minimal point.

Pick any $\alpha \in \RP^1$.  The convex set $\logdom$ has
horizontal and vertical support hyperplanes (by non-entirety of
$F(z,0)$ and $F(0,w)$), and therefore has a support hyperplane
normal to $\alpha$; let $(x,y)$ be a point of intersection of this
support plane with $\logdom$. We have just seen that $\b{z} (\alpha)
:= (e^x , e^y)$ is a minimal point of $\sing$.  If $\b{z}$ is a smooth
point of $\sing$ then $\alpha \in \direc (\b{z})$: either $\b{z}$ is
finitely minimal, in which case Theorem~\ref{th:2d} applies, or it
is toral, in which case the toral version of this theorem from
\citeasnoun{pemantle-wilson;toral} applies.

Assume now that $\b{z}$ is not a smooth point.  By
Lemma~\ref{lem:homogeneous}, $\b{z}$ is a homogeneous point, and
since $d=2$, $\b{z}$ is a multiple point.  Theorem~\ref{th:PW2000a}
then shows that $\alpha \in \direc (\b{z})$ in this case as well.
This finishes the proof. 
\end{proof}

\section{Further details and open questions} \label{ss:open}

The theorems in this and subsequent papers give estimates that are
uniform away from the boundary of the domain in which they are
valid.  In order for all of these to be patched together so as to
give estimates valid now matter how $\b{r} \to \infty$, one must
determine the bandwidth around the boundary for which the boundary
estimates on either side hold.  For instance, suppose $(z,w)$ is a
multiple point of degree 2 and that $\direc (z,w)$ is the set of
slopes between $1/2$ and 2.  It appears that the asymptotic
estimate in \citeasnoun{pemantle-wilson;multiple} holding near the line $\{
s = 2r \}$ can be written so it is valid out to $s = 2r + c
\sqrt{r}$.  If the estimate for the region $s/r > 2 + \ee$ can be
widened so it holds to $s = 2r + c \sqrt{r}$ and a description
given that is valid in the regime $(s-2r)/\sqrt{r} \to c$, then
the estimates will patch together completely.

Another natural question is the universality of the method when
the coefficients have mixed signs.  We conjecture that
Theorem~\ref{th:classify} still holds, in the sense that for every
direction there is point $\b{z} \in \sing$ for which integration near
$\b{z}$ yields correct asymptotics.  What we know is that $\b{z}$ may
no longer be minimal.  For example, if $\numer = 1$ and
$$\denom = (1 - (2/3) w - (1/3) z) (1 + (1/3) w - (2/3) z )$$
then the point $(3/2 , 3/4)$ is not minimal but yields asymptotics
in the diagonal direction; one sees this by integrating along a
deformed torus rather than along $\torus (3/2 , 3/4)$.  In fact we
conjecture that such a deformation always exists, but the topology
seems not transparent enough to yield an easy proof.

The class of algebraic functions is in some ways almost as nice as
the set of rational functions, and nicer than the meromorphic
functions.  For one thing, an algebraic function is determined by
a finite amount of data, and may thus easily be input into a
symbolic math package.
\citeasnoun{gao-richmond;limit-multivariate4} give an analysis of
algebraic and logarithmic singularities, but sometimes the
relevant singularities for algebraic functions are poles.  For
example, in Larsen and Lyons' analysis
\cite{larsen-lyons;coalescing} of merge times for coalescing
particles, they find an algebraic function of the form
$$F(z,w) = {\chi (z,w) \over w - 1 - \sqrt{1-z}}$$
with $\chi$ analytic.  The branch of the square root is chosen so
that at the origin the denominator is 2, not 0.  There is a
branchline at $z=1$, but for all directions in $\RP^1$, there is a
smooth pole on the curve $w = 1 + \sqrt{1-z}$ yielding asymptotics
in the desired direction.  It is natural to ask when this will
happen, and how one can tell effectively.  Some questions of
effectiveness are addressed in \citeasnoun{pemantle-wilson;multiple} and
\citeasnoun{pemantle-wilson;toral},
but there is probably substantial room for improvements on an
algorithmic level.

\noindent{\bf Acknowledgements:} We thank Steve
Wainger for his generous help, including pointing out
equation~(\ref{eq:X}).  Thanks also to Andreas Seeger for his help
with the material in Section~\ref{ss:oscillate}.  We are indebted
to Manuel Lladser for many helpful comments on an earlier draft.

\bibliographystyle{agsm}
\bibliography{asymp}
\end{document}